\documentclass[11pt]{amsart}
\usepackage{amsmath}
\usepackage{amsxtra}
\usepackage{amscd}
\usepackage{amsthm}
\usepackage{amsfonts}
\usepackage{amssymb}
\usepackage{psfig}
\usepackage{eucal}
\newcommand{\psdraw}[3]{\begin{array}{c} \hspace{-1mm}
\raisebox{-4pt}{\psfig{figure=#1.ps,width=#2,height=#3}}
\hspace{-1mm}\end{array}}
\newtheorem{thm} {Theorem}

\newtheorem{cor} {Corollary}

\newtheorem{lem} {Lemma}
\newtheorem{defn} {Definition}

\title[Explicit Enumeration of 321,Hexagon--Avoiding Permutations]{Explicit
Enumeration of 321,Hexagon--Avoiding Permutations}

\author{Zvezdelina Stankova--Frenkel and Julian West}
\address{Zvezdelina Stankova--Frenkel, Dept. of Mathematics and
Computer Science, Mills College, Oakland, CA, {\tt stankova@mills.edu}}
\address{Julian West, Dept. of Mathematics and Statistics, University
of Victoria, Canada, {\tt westj@mala.bc.ca}}

\begin{document}
\thispagestyle{empty}
\setcounter{tocdepth}{1}
\maketitle
\vspace*{-5mm}
\centerline{May 2001}

\begin{abstract} The {\it 321,hexagon--avoiding} ({\it 321--hex})
permutations were introduced and studied by Billey and Warrington in
\cite{Billey} as a class of elements of $S_n$ whose Kazhdan--Lusztig
and Poincar\'{e} polynomials and the singular loci of whose Schubert
varieties have certain fairly simple and explicit descriptions. This
paper provides a 7--term linear recurrence relation leading to an
explicit enumeration of the 321--hex permutations. A complete
description of the corresponding generating tree is obtained as a
by--product of enumeration techniques used in the paper, including
Schensted's 321--subsequences decomposition, a 5--parameter generating
function and the symmetries of the octagonal patterns avoided by the
321--hex permutations.
\end{abstract}

\tableofcontents
\section{Introduction}

We start by describing the {\it 321,hexagon--avoiding} permutations in
two ways; first in the context of pattern--avoidance, which is the
viewpoint we will be exploiting in our enumeration, and then in the
context of reduced expressions, which explains the introduction of the
term 321,hexagon--avoiding (for brevity, {\it 321--hex}). Finally, we
explain briefly the connection with Kazhdan--Lusztig polynomials
in \cite{Billey} and how this motivates the work on the present paper.

\subsection{Pattern--avoidance.}
From the first viewpoint, we consider permutations in $S_n$ as
bijections $w: [n] \rightarrow [n]$, and write them in {one--line
notation} as the image of $w$, $[w_1,w_2,\ldots,w_n]$.  For $n<10$ we
suppress the commas without causing confusion.

\begin{defn}
{\rm Let $v \in S_k$ and $w \in S_n$ for some $k\leq n$. We say that
$w$ {\em contains} $v$ if there is a sequence $1 \leq i_1 < \cdots <
i_k \leq n$ such that the sequences
$w^{\prime}=[w_{i_1},w_{i_2},\ldots, w_{i_k}]$ and
$[v_1,v_2,\ldots,v_k]$ obey the same pairwise relations, i.e.
$w_{i_j}<w_{i_m}$ exactly when $v_j<v_m$.  In such a case, we write
$w^{\prime}\sim v$. If $w$ does not contain $v$ then we say that $w$
{\em avoids} $v$. We denote by $S_n(v)$ the set of all $v$--avoiding
permutations of length $n$.}
\end {defn}

\noindent For example, the permutation $\omega=(52687431)$ avoids
$(2413)$ but does not avoid $(3142)$ because of its subsequence
$(5283)$.  For a classification of forbidden subsequences
up to length 7, we direct the reader to \cite{BW,BWX,St1,St2,St3,We1,We2}.

\begin{defn}
{\rm The {\it 321--hex} permutations are those permutations which
simultaneously avoid each of the following 
five patterns: 
\[[321],\,P_1=[46718235],\,P_2=[46781235],\,P_3=[56718234],\,P_4=[56781234].\]}
\end{defn}
\noindent We denote by $\mathcal{P}$ the set of the four length--8
({\it octagonal}) permutations $P_1,P_2,P_3,P_4$.  In order to make
sense of the above definition, consider the following equivalent, but
perhaps more insightful, reformulation in terms of matrices.

\begin{defn}
{\rm Let $w \in S_n$. The {\it permutation matrix} $M_w$ is the
$n\times n$ matrix having a $1$ in position $(i,w_i)$ for $1\leq i
\leq n$, and $0$ elsewhere. (To keep the resemblance with the
``shape'' of $w$, we coordinatize $M_w$ from the bottom left corner.)
Given two permutation matrices $M$ and $N$, we say that $M$ {\it
avoids} $N$ if no submatrix of $M$ is identical to $N$.}
\end{defn}

\begin{figure}[h]
$$\psdraw{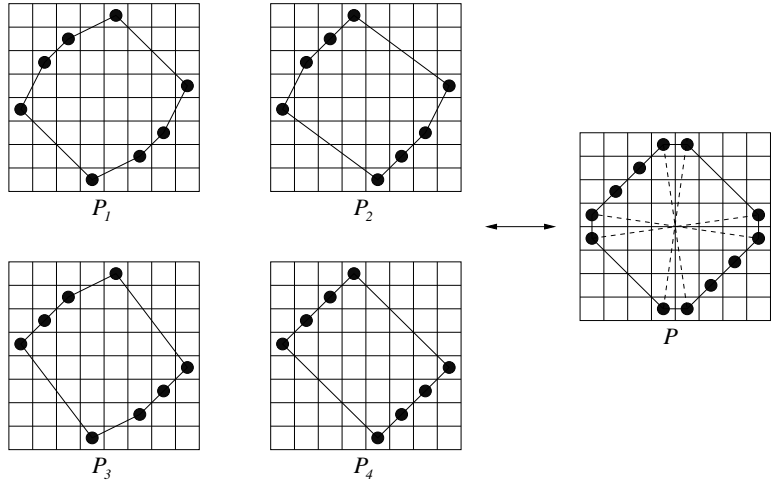}{4in}{2.5in}$$
\caption{Octagonal Patterns}
\label{octagonal}
\end{figure}

A permutation matrix $M$ of size $n$ is simply a transversal of an
$n\times n$ matrix. Clearly, $w \in S_n$ contains $v \in S_k$ if and
only if $M_w$ contains $M_v$ as a submatrix. Under this reformulation,
Fig.~\ref{octagonal} presents the four octagonal patterns
$\{P_i\}$ which must be avoided by the 321--hex
permutations. The fifth pattern $P$ does {\it not} come from a
permutation because it is not a transversal. Yet, $P$ is a union of
the four previous permutation patterns, and it can be easily checked
that a permutation matrix $M_w$ avoids all $P_i$'s if and only if no
$8\times 8$ permutation submatrix of $M_w$ can be completely
``covered'' by $P$. Thus, by abuse of notation, we can say that the
321--hex permutations are defined as the permutations avoiding both
$[321]$ and $P$.

\smallskip
Our enumeration makes use of this last interpretation, exploiting the
symmetries in the set $\mathcal{P}$ of octagonal patterns
plus a convenient structural representation of all 321--avoiding
permutations.  The close relations between the four octagonal
permutations are even more clearly revealed from a group--theoretical
viewpoint when examining their reduced expressions. Below, we briefly
review the construction of the {\it heaps} of 321--avoiding
permutations and their relation to the octagonal patterns $P_i$; for
more details, see \cite{Billey}. The remainder of the introduction,
except for its conclusion, can be skipped by the reader who is
interested only in the pattern--avoidance interpretation.

\subsection{Heaps of 321--permutations.}
The permutation group $S_n$ can be regarded as generated by the set of
the {\em adjacent transpositions} $\{s_i\}_{i=1}^{n-1}$, where
$s_i=(i,i+1)$ in {\it cyclic} notation. In this presentation, the
generators $s_i$ and $s_j$ commute if $|i-j|>1$; else $s_is_{i+1}s_i =
s_{i+1}s_is_{i+1}$.  An {\em expression} is any product of generators.
A {\em reduced} expression $\langle w \rangle$ for $w\in S_n$ is a
shortest--possible expression yielding $w$. (It is well--known that the
number of generators in $\langle w \rangle$ equals the number of
inversions in $w$.)  For example, the octagonal pattern
$P_1=[46718235]$ has a reduced expression
\begin{center}
$\langle  P_1 \rangle =s_3s_2s_1s_5s_4s_3s_2s_6s_5s_4s_3s_7s_6s_5.$
\end{center}
Reduced expressions for the other avoided octagonal patterns are
$\langle P_2 \rangle=\langle P_1 \rangle \cdot s_4$, $\langle P_3
\rangle=s_4\cdot\langle P_1 \rangle$ and $\langle P_4
\rangle=s_4\cdot\langle P_1 \rangle\cdot s_4$. These expressions can
be easily verified by considering the dashed lines in $P$ in
Fig.~\ref{octagonal}. The study of reduced expressions is a major
subject of the representation theory of $S_n$.

\smallskip
After work of Viennot \cite{Vi}, the $321$--avoiding permutations $w$
can be represented by special ranked posets called {\em heaps}.  The
elements of $\text{Heap}(w)$ are identified with the transpositions
$\{s_{i_j}\}$ in a fixed reduced expression $\langle w\rangle$ for
$w$. By Billey--Jockusch--Stanley~\cite{BJS}, the $321$--avoiding
permutations are those in which no reduced expression contains a
substring of the form $s_is_{i\pm 1}s_i$; and by Tits \cite{Ti}, all
reduced expressions for a $321$--avoiding permutation $w$ are
equivalent up to moves $s_is_j \rightarrow s_js_i$ for
$|i-j|>1$. Thus, the set of elements in $\text{Heap}(w)$ is
independent of the choice of $\langle w\rangle$.

\begin{figure}[h]
$$\psdraw{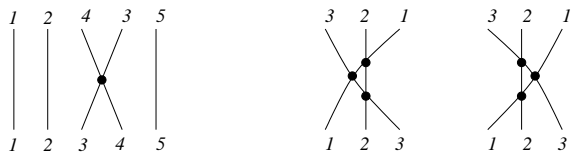}{3in}{0.75in}$$
\caption{String diagrams for $s_3\in S_5$, $s_2s_1s_2\in S_3$ and
$s_1s_2s_1\in S_3$}
\label{stringdiagram}
\end{figure}

We now describe the {\it rank function} of $\text{ Heap}(w)$, along
with a Hasse diagram for the poset by embedding its elements in the
integer lattice. One way to define and visualize this embedding is via
{\it string diagrams}.  To form a string diagram of $w\in S_n$, write
the row of numbers $[w_1,w_2,...,w_n]$ above the row $[1,2,...,n]$,
thus mimicking the two--line notation for a permutation.  Connect each
number $i$ on the bottom line to the corresponding number $i$ on the
top line, drawing a "string" which may change direction, but which at
all times is running either due north, northwest or northeast.
Strings may cross and recross, but do not run over top of one another,
nor do they stray beyond the rectangular bounds formed by the
two rows of numbers.  For example, $s_i\in S_n$ can be
realized as the crossing of two strings, as shown in
Fig.~\ref{stringdiagram}a.  In general, two or more adjacent
transpositions can be applied simultaneously, provided they commute,
i.e. $s_i$ and $s_j$ can occur on the same horizontal level of a
string diagram unless $|i-j| = 1$.

Thus, the crossings of a string diagram can be identified with the
$s_i$'s, labelled by the column (from $1$ to $n-1$) in which they
occur, and can be seen to form a poset in the obvious
way.  The linear extensions of this partial order are the expressions
for $w$ defined above. If the string diagram has the smallest possible
number of crossings, then it is is {\em minimal} and the linear
extension is a reduced expression for $w$.

A {\em short braid} is a configuration obtained by applying the
non--commuting transpositions $s_i$ and $s_{i+1}$ in the following
order: $s_is_{i+1}s_i$ or $s_{i+1}s_is_{i+1}$
(cf.~Fig.~\ref{stringdiagram}b--\ref{stringdiagram}c).  The string
diagram for a braid shows crossings at three of the four points of a
small diamond, omitting either the eastern or the western point.  As
mentioned earlier, a permutation is 321--avoiding exactly when its
minimal string diagrams avoid such configurations; such permutations
are therefore also called {\em short--braid--avoiding} in the
literature of Coxeter groups. It is now clear how to embed canonically
the poset of string crossings in a minimal string diagram of $w$ into
the integer lattice.  The resulting Hasse diagram is the heap of $w$,
$\text{Heap}(w)$. It is independent of the choice of $\langle
w\rangle$ as long as $w$ is 321--avoiding.
 
\begin{figure}[h]
$$\psdraw{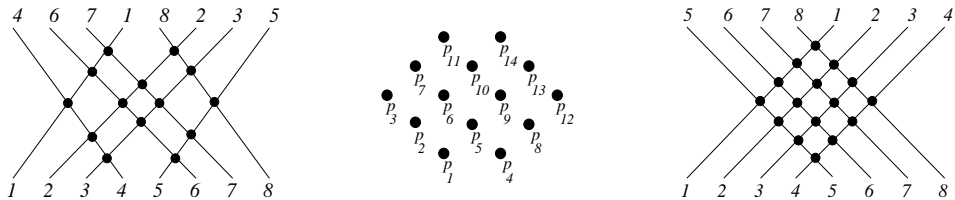}{5in}{1in}$$
\caption{String diagrams for $P_1$ and $P_4$, and 
Heap$\langle P_1 \rangle$ where ${\langle P_1
\rangle=\prod_{i=1}^{14}p_i}$}
\label{heap}
\end{figure}

\smallskip
The heap of the special reduced expression $\langle P_1 \rangle$ for
the octagonal pattern $P_1$ resembles a hexagon: ${\rm Heap}\langle
P_1 \rangle$ has horizontal and vertical symmetries, with respectively
2,3,4,3,2 lattice points on its five ranks (see
Fig.~\ref{heap}a--\ref{heap}b). The string diagram of $P_4$ features
one extra point on top and one extra point on bottom, corresponding to
the crossings of the strings $1$ and $8$, and $4$ and $5$ (see
Fig.~\ref{heap}c). The string diagrams for $P_2$ and $P_3$ have either
the top or the bottom extra crossing. In all cases,
$\text{Heap}\langle P_i \rangle$ contains the hexagonal
$\text{Heap}\langle P_1 \rangle$.

Furthermore, it can be shown that the $(321,\mathcal{P})$--avoiding
permutations are exactly those 321--avoiding permutations whose
minimal string diagrams avoid the hexagonal string diagram of
$P_1$. This justifies the descriptive term {\it 321,hexagon--avoiding}
and yields the following alternative description (cf.~\cite{Billey}):

\smallskip
\noindent{\bf Definition $\bf 2^{\prime}$.} A {\it 321--hex
permutation} is a permutation whose reduced expressions do not contain
a substring of the form $s_js_{j\pm1}s_j$ for all $j\geq 1$, and do
not contain any of the translates ${\prod_{k=1}^{14}s_{i_k+j}}$ of
$\text{Heap}\langle P_1 \rangle$ for all $j\geq 0$.


\subsection{Kazhdan--Lusztig polynomials.}
We now turn our attention to the Kazhdan--Lusztig polynomials.  For
definitions and a detailed introduction, we direct the reader to
\cite{Hum},\cite{KL}. In short, the {\it Hecke algebra} $\mathcal{H}$
of a finite Weyl group $W$ (such as $S_n$) is an algebra over 
$\mathbb{Q}(\sqrt{q})$, with basis $\{T_w\}_{w\in W}$, relations 
for all generators $s\in W$:
\[\left|\begin{array}{rcl}
T_sT_w&=&T_{sw}\,\,\,\,\text{if}\,\,\,\,l(sw)>l(w),\\
T_s^2&=& (q-1)T_s+qT_1,
\end{array}\right.\]
and well--defined inverses due to the presence of $q^{-1}$. Moreover,
the involution of $\mathbb{Q}(\sqrt{q})$ sending $\sqrt{q}\mapsto
1/\sqrt{q}$ extends to an involution $\iota$ of $\mathcal{H}$. The
{\it Kazhdan--Lusztig polynomials} arise in search of a new basis
$\{C_w^{\prime}\}$ for $\mathcal{H}$ of $\iota$--invariant elements
which are linear combinations of ``lower--terms'' $T_x$ for $x\leq w$
under the Bruhat--Chevalley order: $x\leq w$ if every reduced
expression for $w$ contains a subexpression for $x$.

\begin{thm}[Kazhdan--Lusztig] For any $w\in W$, $\mathcal{H}$ has
a unique $\iota$--invariant element
$C_w^{\prime}=q^{-l(w)/2}\sum_{x\leq w}P_{x,w}T_x$, where the degrees
of the polynomials $P_{x,w}(q)\in \mathbb{Z}[q]$ are at most
$\frac{1}{2}(l(w)-l(x)-1)$ if $x<w$, and $P_{w,w}=1$, $P_{x,w}=0$ if
$x\not \leq w$.
\end{thm}

The Kazhdan--Lusztig polynomials $P_{x,w}$ are of fundamental
importance in Lie Theory. It has been proven that their coefficients
are non--negative for Weyl groups (see Kazhdan--Lusztig \cite{KL2},
Brylinski--Kashiwara \cite{BK}), but this question remains open for
arbitrary Coxeter groups.  Neither the degrees nor the coefficients of
$P_{x,w}$ are readily computable; indeed, only partial results have
been derived in certain cases (see \cite{Billey} for detailed
references.)

\subsection{321--hex permutations in Kazhdan--Lusztig polynomials.}
In \cite{Billey}, Billey and Warrington introduce the 321--hex
permutations $w\in W=S_n$, derive simple combinatorial formulas for
their Kazhdan--Lusztig and Poincar\'{e} polynomials, and a simple
method for determining the singular loci of their Schubert varieties
$X_w$. Concretely, let $\langle w \rangle=s_{i_1}\cdots s_{i_r}$ be a
reduced expression for $w$, and let $d(\sigma)$ denote the {\it defect
statistic} of a {\it mask} $\sigma$ on $\langle w \rangle$ whose
product is a fixed $x\in S_n$. (For definitions of these and other
related concepts, see \cite{Billey,BS}.)

\begin{thm}[Billey--Warrington] A permutation $w\in S_n$ is
321--hexagon avoiding if and only if one of the following
equivalent conditions is satisfied:
\begin{itemize}
\item[(a)] The Kazhdan--Lusztig polynomial 
${P_{x,w}=\sum_{\sigma}q^{d(\sigma)}}$ for $x\leq w$.

\item[(b)] The Poincar\'{e} polynomial for the full intersection
cohomology group of $X_w$ equals $(1+~q)^{l(w)}$.

\item[(c)] The Kazhdan--Lusztig basis element $C^{\prime}_w=
C^{\prime}_{s_{i_1}}\cdots C^{\prime}_{s_{i_r}}.$

\item[(d)] The Bott--Samelson resolution $Y$ of $X_w$ is small.

\item[(e)] $\text{IH}_*(X_w)\simeq H_*(Y)$.
\end{itemize}
\label{BW-theorem}
\end{thm}

\subsection{Enumeration of 321--hex permutations.} As 
Theorem~\ref{BW-theorem} shows, the 321--hex permutations have
properties that make certain algebraic and combinatorial computations
easier to perform compared to arbitrary permutations. Another
interesting aspect of the 321--hex permutations was examined earlier
in the introduction: they are one of the few families that are known
to be describable both in terms of pattern--avoidance and
heap--avoidance. Yet, until now, there was no known recursive, exact
or other closed form for the number of 321--hex permutations. For
instance, it would be useful to know to how many permutations
Theorem~\ref{BW-theorem} applies; how the number of 321--hex
permutations changes asymptotically, and how it compares to sizes of
other well--known sets of permutations, such as $S_n(321)$, which is
enumerated by the Catalan numbers.

\smallskip
Answers to all these questions are obtained in the present paper,
where we find a 7--term linear recursive relation and derive from it an
explicit exact formula. With this formula at hand, answering any
enumeration questions about the 321--hex permutations becomes a matter
of simple observation and calculation.

\begin{defn} Let $\mathcal{H}_n$ denote the set of all 321--hex
permutations in $S_n$, and let $\alpha_n=|\mathcal{H}_n|$ be the
number of such permutations.
\end{defn}

\begin{thm} The sequence $\alpha_n$ satisfies the following
recursive relation for all $n\geq 6$:
\begin{center}
$\alpha_{n}=6\alpha_{n-1}-11\alpha_{n-2}+
9\alpha_{n-3}-4\alpha_{n-4}-4\alpha_{n-5}+\alpha_{n-6}.$
\end{center}
\label{main-result}
\end{thm}
\noindent Of the six roots of the corresponding characteristic
polynomial, four are real: $R_i$ for $i=1,2,3,4$, and two are complex
conjugates: $R_5=\overline{R_6}$. This implies the same description of
the six coefficients below: $c_i\in\mathbb{R}$ for $i=1,2,3,4$, and
$c_5=\overline{c_6}\in \mathbb{C}$.

\begin{cor} The number of the 321--hexagon
avoiding permutations of length $n$ equals
\[c_1R_1^n+c_2R_2^n+c_3R_3^n+c_4R_4^n+c_5R_5^n+
\overline{c_5{R_5}^n},\]
where the roots and coefficients are rounded off below 
to 5 digits after the decimal point:
\[\begin{array}{lll}
R_1\approx -0.49890          & &c_1\approx  0.00164\\
R_2\approx\phantom{-}0.21989& &c_2\approx  0.13776\\
R_3\approx\phantom{-}1.95627& &c_3\approx  0.57156\\
R_4\approx\phantom{-}3.43526& &c_4\approx  0.24149\\
R_5\approx\phantom{-}0.44375-1.07682i& &c_5\approx  0.02378+0.00080i\\
\end{array}\]

\end{cor}
\noindent For further discussion of this and other results, we refer
the reader to Sections~\ref{results}--\ref{extensions}. The proof of
Theorem~\ref{main-result} follows several steps. First, we describe
the nodes in the generating tree of $\mathcal{H}_n$ by using
Schensted's algorithm for 321--avoiding permutations and by
introducing 5 parameters for the generating function
$h_n(x,k,l,m)$. We next observe that this function depends on fewer
parameters, yielding therefore relatively few distinct values. We
organize these values in five sequences,
$\alpha_n,\beta_n,\gamma_n,\delta_n$ and $\epsilon_n$. Using the
intrinsic symmetries of the set $\mathcal{P}$ of octagonal patterns,
we deduce recursive relations expressing each sequence in terms of
$\alpha_n$. The latter turns out to be the number we are looking for,
$|\mathcal{H}_n|$. Finally, putting together all information about the
generating function and the five sequences results in the desired
formula.

\section{The generating tree}
\label{generating}

\subsection{What is a generating tree?} 
We turn now to the development of a recurrence for the 321--hex
permutations.  A standard tool in the enumeration of restricted
permutations is the {\em generating tree} $T$ introduced in
\cite{CGHK}.  Begin with an infinite tree whose nodes on level $n$ are
identified with the permutations in $S_n$. The node $w$ is a child of
$\widehat{w} = [w_1,\ldots,w_{j-1},w_{j+1},\ldots,w_n]$ where the
omitted value is $w_j=n$.  Looking at this from the point of view of
the parent, we can form all the children of $w\in S_{n}$ by {\em
inserting} the element $n+1$ into each of the $n+1$ {\em sites} of
$w$.

\smallskip
Then, for a given set $\Sigma$ of forbidden subsequences, prune the
tree by deleting all nodes containing any of the forbidden
subsequences.  What remains is still connected because if $w$ does not
contain any forbidden subsequence then clearly $\widehat{w}$ does not
either.  For any node on level $n$ of this pruned tree $T(\Sigma)$, we
call a site in the corresponding permutation {\em active} if
inserting $n+1$ at that site yields a node of the tree; conversely an
{\em inactive} site is one where the insertion of $n+1$ creates one or
more of the forbidden subsequences in $\Sigma$.

\smallskip
To give a complete description of a generating tree, we need to
associate to each node an appropriate label, and then describe a {\em
succession rule} for deriving the labels attached to the set of
children of each node.  For instance, we might characterize the
original tree $T$ generating {\em all} permutations as having a root
labelled $(2)$ and a succession rule $(n) \rightarrow (n+1)^n$.  In
this instance, the label can be interpreted as revealing directly how
many children each node has in the generating tree.

\smallskip
In the particular instance of the 321--hex permutations, it will turn
out that we need a label containing four integers.  Although this is
more complicated than the single integer of our motivating example, it
is nevertheless a major progress to reduce the amount of information
recorded at a node from a full permutation to a label of any bounded
size.  In particular, it is possible to apply the succession--rule
recursively to determine the entire downward structure of any node
given only its label, regardless of whether that node is on a level
corresponding to permutations on four symbols, or four thousand.

\subsection{Schensted algorithm for $S_n(321)$ and active regions.} 
Let $w$ be any 321--hexagon--avoiding permutation on $n$ symbols.  We
divide the elements $w_1,w_2,\ldots,w_n$ into two categories: the set
of right--to--left minima (including $w_n$), and the rest.  Adapting the
terminology of Schensted \cite{Schensted}, we refer to these as the
{\em first basic subsequence} and {\em second basic subsequence} of
$w$, and denote them by $\mathcal{B}_1(w)$ and $\mathcal{B}_2(w)$,
respectively.  For example, when $w=P_1$, $\mathcal{B}_1(w)=[4,6,7,8]$
and $\mathcal{B}_2(w)=[1,2,3,5]$.  Note that both $\mathcal{B}_1(w)$
and $\mathcal{B}_2(w)$ decrease from right to left, the former -- by
construction, and the latter -- because $w\in S_n(321)$.  Now let
$K=w_{i_1}$, $L=w_{i_2}$ and $M=w_{i_3}$ be the three largest
(i.e. rightmost) elements in $\mathcal{B}_2(w)$ ($K<L<M$), with each
one set to zero if the corresponding element does not exist.

\smallskip
Since $M$ is the largest (i.e. rightmost) element in
$\mathcal{B}_2(w)$, it follows that every element to the right of $M$
belongs to $\mathcal{B}_1(w)$.  There are $x:=n-i_3$ elements in this
region, which we will call the {\em active region}.  If $M=0$ because
$\mathcal{B}_2(w)$ is empty, then we consider the entire permutation
to be the active region.  Let $k$ be the number of elements in the
active region which are larger than $K$; as the elements in the active
region decrease from right to left, these $k$ elements are
$w_{n-k+1},\ldots,w_{n}$.  Similarly, $l$ be the number of elements in
the active region which are larger than $L$, and $m$ be the number of
such elements larger than $M$ (cf. Fig.~\ref{deletion1}.)

Now assign the label $(x,k,l,m)$ to $w$.  By construction, $x\geq k
\geq l\geq m$. Let the {\it $X$--elements} of $w$ be those that are
counted in $x$, and similarly define the $K$--, $L$--, and
$M$--elements. Further, let the {\it $X\backslash K$--elements} of $w$
be the set of $X$--elements minus the set of $K$--elements, and
similarly define the $K\backslash L$-- and $L\backslash M$--elements.
For example, the $X$--elements in $P_1$ are $2,3,5$, and $k=l=m=0$.
\subsection{The succession rule for 321--hex permutations.} 
All children of $w$ must avoid the subsequence $[321]$.  This
restriction by itself renders inactive all sites to the left of $M$,
but none of the sites in the active region.  Therefore when
considering the children of $w$ and their labels, we need only
consider insertions taking place in the active region to the right of
$M$.  Thus, if we insert $n+1$ into a site in the active region with
$j$ elements to its right, it is easy to verify that the resulting
permutation will have label:
\begin{equation}
\left\{\begin{array}{cl}
(x+1,k+1,l+1,m+1)& \text{if}\,\, j=0\\
(j,j,j,0)        & \text{if}\,\, 0 < j \leq m\\
(j,j,m,0)        & \text{if}\,\, m < j \leq l\\e
(j,l,m,0)        & \text{if}\,\, l < j \leq x\\
\end{array}\right\}
\label{labels1}
\end{equation}
Furthermore, the number of $w$'s children can be computed as follows.
Set \[T:=\min(k+2,\max(k+1,l+2)).\] Then the node $w$ has $S+1$
children, corresponding to the $S+1$ rightmost insertion sites, where
\begin{equation}
\left|\begin{array}{ll}
S=T&\text{if}\,\,T\leq x-2,\\
S=x&\text{if}\,\,T> x-2.\\
\end{array}\right.
\label{children}
\end{equation}
This allows a more compact and complete succession rule for the labels
of all $S+1$ children of $w$:
\begin{equation}
(x,k,l,m)\,\,\mapsto\,\,\left\{\begin{array}{l}
                          (x+1,k+1,l+1,m+1),\\
                          \big(i,\min(i,l),\min(i,m),0\big)\,\,
                          \text{for}\,\,i=1,...,S.
\end{array}\right\}
\label{labels2}
\end{equation}
We shall not use directly any of the above formulas
(\ref{labels1})--(\ref{labels2}) in our calculations, so we leave
their verification to the reader, who will find this easier after
mastering the material in the rest of the paper. The importance of the
above discussion is that it completely describes the structure of the
generating tree $T(321,\mathcal{P})$, and hence explains in principle
why the enumeration in this paper works. Why the resulting final
formula for the 321--hex permutations is so simple -- a {\it linear}
recursive relation with {\it constant} coefficients -- is a completely
different matter and can be explained only by the structure of the
forbidden set $\mathcal{P}$ of four octagonal patterns, as we shall
see later.

\subsection{The generating function for $T(321,\mathcal{P})$.}
Let $h_n(x,k,l,m)$ and $\mathcal{H}_n(x,k,l,m)$ be the number,
resp. the set, of 321--hex permutations in $S_n$ labelled
$(x,k,l,m)$. According to this definition, $h_n(0,0,0,0)$ is not
defined since no permutation is labelled $(0,0,0,0)$. For convenience,
denote by $h_n(0,0,0,0)$ and $\mathcal{H}_n(0,0,0,0)$ the number and
the set of all 321--hex permutations in $S_n$ which end in their
largest element: $w_n=n$. It is worth noting that
$\mathcal{H}_n(x,x,x,m)$ for $x>m$ corresponds to permutations $w$ in
which either $L$ is smaller than the final ``tail'' of $w$, or $L$
does not exist and $\mathcal{B}_2(w)=\{M\}$.

\medskip
Naively, there are $\mathcal{O}(n^5)$ enumerations to be performed for
$h_n(x,k,l,m)$. However, many have the same answers; indeed, only $5n$
different values of $h_{n+1}(x,k,l,m)$ appear:

\[\begin{array}{|c|r|r|r|r|r|r|r|r|r|r|r|r|r|r|r|}\hline
  n & 1 & 2 & 3  & 4 & 5 & 6 & 7 & 8 & 9 & 10 & 11 & 12\\ \hline\hline
    \alpha_n  &1&2&5&14&42&132&429&1426&4806&16329&55740&190787\\
    \beta_n   &0&0&1&4&14&48&165&568&1954&6717&23082&79307\\
    \gamma_n  &0&0&0&1&5&20&75&271&957&3337&11559&39896\\
    \delta_n  &0&0&0&0&1&6&25&93&333&1172&4083&14137\\
    \epsilon_n&0&0&0&0&0&1&5&19&68&240&839&2911\\\hline
\end{array}\]

\begin{lem} The operation of deleting all $K$--elements 
in a 321--hex permutation provides a bijection
$d_K:\mathcal{H}_n(x,k,l,m)\stackrel{\sim}{\rightarrow}
\mathcal{H}_{n-k}(x-k,0,0,0)$. Hence, for all $n,x,k,l,m$ we
have $h_n(x,k,l,m)=h_{n-k}(x-k,0,0,0)$.
\label{deletion1}
\end{lem}

\noindent{\sc Proof:} Let $w\in S_n(321)$. The $K$--elements of $w$
lie in $\mathcal{B}_1(w)$ and are part of the final increasing
``tail'' of $w$ (since they are to the right of $M$). Thus, if a
$K$--element $w_i$ were part of an octagonal pattern $P_j$ in $w$,
then $w_i$ would lie in $\mathcal{B}_1(P_j)$. But then $P_j$, and
hence $w$, would contain at least 3 elements {\it larger than} and
{\it to the left of} $w_i$: this contradicts the definition of a
$K$--element, for which only $L$ and $M$ are {\it larger than} and
{\it to the left of} it.

\begin{figure}[h]
$$\psdraw{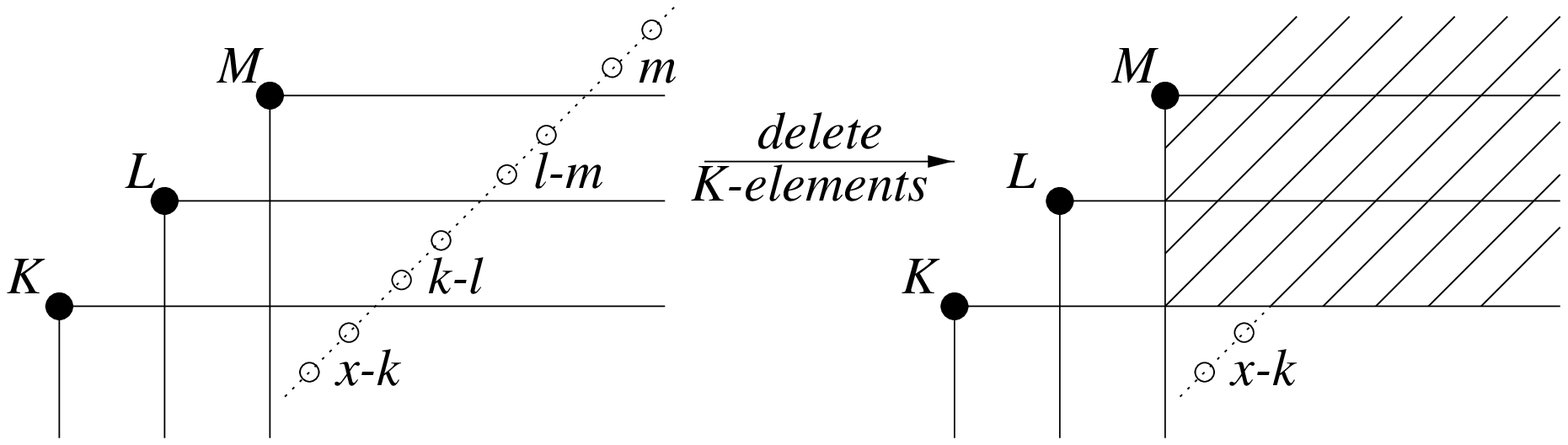}{4in}{1.2in}$$
\caption{Lemma \ref{deletion1}}
\end{figure}

\smallskip
The above discussion shows that $K$--elements cannot participate in
octagonal patterns in $w$, and hence they can be deleted without
losing any relevant 321--hex information about $w$.  (Of course, we
have  to rescale down appropriately $M$ and $L$ of $w$ to arrive at a
permutation of smaller size.) The resulting map
\[d_K:\mathcal{H}_n(x,k,l,m){\rightarrow} \mathcal{H}_{n-k}(x-k,0,0,0)\]
is bijective: to obtain $w\in \mathcal{H}_n(x,k,l,m)$ from its image
$\widetilde{w}=d_K(w)\in\mathcal{H}_{n-k}(x-k,0,0,0)$, insert the
necessary number of $M,L,K$--elements into $\widetilde{w}$ and
increase appropriately $L$ and $M$ to fit their definitions in $w$.
This procedure works because we can identify $M$ as the largest
element in $\widetilde{w}$ (after the deletion of all $K$--elements in
$w$), and then identify $L$ as the second largest element of
$\widetilde{w}$, which will be necessarily to the left of $M$. Then
insertion of the appropriate $M,L,K$--elements at the end of $w$
requires rescaling $M$ and $L$ only (not $K$), and hence transforms
the 321--Schensted decomposition of $\widetilde{w}$ into that of $w$:
whether $L$ was in $\mathcal{B}_1(\widetilde{w})$ or in
$\mathcal{B}_2(\widetilde{w})$ does not prevent $L$ from becoming an
element of $\mathcal{B}_1(w)$ after applying $d_K^{-1}$.  \qed

\subsection{Octagonal conditions.}
Consider a family $\mathcal{F}$ of permutations in
$\mathcal{H}_n$ which are described by certain configuration
conditions imposed on their basic subsequence decomposition
$\mathcal{B}_1\sqcup \mathcal{B}_2$. For example,
$\mathcal{H}_n(0,0,0,0)$ is such a family defined by the condition
$w_n=n\in\mathcal{B}_1(w)$. Let $a$ be an element of the permutations
in $\mathcal F$ which is identified uniquely in each $w\in\mathcal{F}$
by the configuration description of $\mathcal F$.

\begin{defn} {\rm We say that $a$ is {\it pattern--free} if its
deletion in each $w\in\mathcal{F}$ (and appropriate rescaling of $w$)
results in a numerically equivalent family $\widetilde{\mathcal{F}}$
of permutations in $\mathcal{H}_n$; i.e. $d_a:{\mathcal
F}{\rightarrow} \widetilde{\mathcal{F}}$ with
$|\mathcal{F}|=|\widetilde{\mathcal{F}}|$. For $w\in \mathcal{F}$,
denote by $\widetilde{w}$ the image
$d_a(w)\in\widetilde{\mathcal{F}}$.}
\end{defn}

\noindent For example, in $\mathcal{H}_n(0,0,0,0)$, the largest
element $a=n$ is identified by being in the last position in each $w$,
and clearly it is pattern--free (no octagonal pattern $P_i$ has $8$ in
its last position):
\[d_n:\mathcal{H}_n(0,0,0,0)\stackrel{\sim}{\rightarrow}\mathcal{H}_{n-1}.\] 
Establishing pattern--free elements and identifying the image set
$\widetilde{\mathcal{F}}$ is the basis of the enumeration of
$\mathcal{H}_n$. The following technical lemma summarizes the
pattern--free situations which will be used later in the proof of
Theorem~\ref{main-result}.

\begin{lem} Let $\mathcal{F}$ be a family in $\mathcal{H}_n$,
$a_1,...,a_{x-k}$ be  the $X\backslash K$--elements of
$w\in\mathcal{F}$, and $H$ be the fourth largest element in
$\mathcal{B}_2(w)$.  The set $\mathcal P$ imposes the following
octagonal conditions on $\mathcal{F}$:

\begin{itemize}
\item [(P1)] All $K$--elements are pattern--free.

\item [(P2)] If $w_n=n$ in $\mathcal{F}$, or $w_n=n-1$ in
$\mathcal{F}$, then $w_n$ is pattern--free.

\item [(P3)] If $M=n$ and $x\leq 2$ in $\mathcal{F}$, then $M$ is
pattern--free.
\end{itemize}

Assume now that $x\geq 3$ and  $k=l=m=0$ in $\mathcal{F}$.

\begin{itemize}
\item [(P4)] If $H<a_2$ in $\mathcal{F}$ ($H$ may not exist in $w$,)
then $M,a_2,...,a_x$ are pattern--free.

\item [(P5)] If $a_2<H<a_3$ in $\mathcal{F}$, then $L,M,a_3,...,a_x$
are pattern--free.

\item [(P6)] If $H>a_3$ (forcing $x=3$) in $\mathcal{F}$, then $M$ is
pattern--free.
\end{itemize}
\end{lem}

\begin{figure}[h]
$$\psdraw{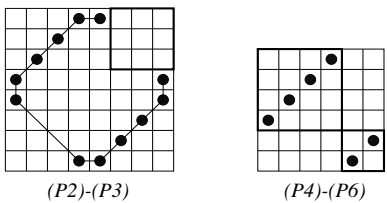}{2in}{1in}$$
\caption{Octagonal Conditions}
\label{octagonal-conditions}
\end{figure}

\noindent{\sc Proof:} (P1) follows from the proof of
Lemma~\ref{deletion1}.  (P2) says that if all $w\in\mathcal{F}$ end in
their largest element $n$ (or in their second largest element $n-1$),
then this last element is pattern--free (cf. Fig.~\ref{alpha});
(P3) says that if $M$ happens to be the largest element in all
$w\in\mathcal{F}$ and is followed by at most two $X$--elements, then
$M$ is pattern--free (cf. Fig.~\ref{beta}--\ref{gamma}.)  Both (P2)
and (P3) follow from the facts:

\begin{itemize}
\item[$\bullet$] Each octagonal pattern $P_i$ has an empty $3\times 3$
upper--right corner (see Fig.~\ref{octagonal-conditions}a.)

\item[$\bullet$] For (P2), reinserting $n$ or $n-1$ in the last
position in $\widetilde{w}\in\mathcal{H}_{n-1}$ does not create
$321$--patterns. Similarly, for (P3), reinserting $M=n$ into
$\widetilde{w}\in\mathcal{H}_{n-1}$ as the first (rightmost) element
in $\mathcal{B}_2(w)$ does not create $321$--patterns (since the tail
of ${w}$ after $M$ is increasing as part of
$\mathcal{B}_1(\widetilde{w})$.)
\end{itemize}

\smallskip
\noindent For (P4), cf.~Case~1~in~Fig.~\ref{epsilon2}. Start by
deleting the crosses of the elements in the bottom row and in the
rightmost column of each $P_i$; this leaves the permutation matrix
$M[345612]$, which decomposes into a $4\times 4$ block $I_4$ in the
upper left corner and the fixed $2\times 2$ block $I_2$ in lower left
corner (corresponding to the original $p_6=2$ and $p_7=3$, see
Fig.~\ref{octagonal-conditions}b.)

\begin{itemize}
\item[$\bullet$] The fact that $I_4$ is to the left of and higher than
$I_2$ proves that $a_2,...,a_x$ are pattern--free: for them, only
$M,L,K$ can fit in such an $I_4$. In detail, if one of $a_2,...,a_x$
participates in an octagonal pattern $P_i$, then without loss of
generality we may assume that $a_x$ is the rightmost element of
$P_i$. This forces $M,L,K$ to participate in $P_i$ too as the only
elements larger than $a_x$. Now $M$, being the largest element of
$P_i$, forces at least three elements after it to participate in $P_i$;
one is $a_x$, one could be $a_1$, hence the third one must be among
$a_2,...,a_{x-1}$. But none of these elements can fit into an
$I_2$ as in Fig.~\ref{octagonal-conditions} because only $L,M,N$ are
larger and to the left of them. Further, reinsertion of $a_2,...,a_x$
into $\widetilde{w}$ as the largest elements of $\mathcal{B}_1(w)$
cannot create $321$--patterns.  Thus, $a_2,a_3,...,a_x$ are indeed
pattern--free.

\item[$\bullet$] After deleting $a_2,...,a_x$, the largest element $M$
lands in second to last position in $\widetilde{w}$, and by (P3), it
is pattern--free.
\end{itemize}

\smallskip
\noindent For (P5), cf.~Case~2 in Fig.~\ref{epsilon2}. The reasoning here is
similar to the case for (P4).  

\begin{itemize}
\item[$\bullet$] First note that there can be no element $a_0$ of $w$
between $L$ and $M$, or else $[H,K,L,a_0,M,a_1,a_2,a_3]\sim
P_1$. 

\item[$\bullet$] If one of $a_3,...,a_x$ participates in a pattern
$P_i$, this forces $M,L,K$ to participate too, and in order not to run
into contradiction with the $I_4\times I_2$ argument above, we must
assume that $a_1$ and $a_2$ are also in $P_i$, and only one among
$a_3,...,a_x$ is in $P_i$. Then the ''1'' in $P_i$ will have to be
between the two largest elements $L$ and $M$, which was ruled
earlier. Thus, $a_3,...,a_x$ are pattern--free.

\item[$\bullet$] Deletion of $a_3,...,a_x$ leaves $M\in
\mathcal{B}_2(\widetilde{w})$ in third to last position, so by (P3),
$M$ is pattern--free. But $L$ is the largest element in
$\widetilde{w}=d_M(w)$, with only $a_1$ and $a_2$ after it, so by
(P3) again, $L$ is pattern--free.
\end{itemize}

\noindent For (P6), cf.  Case~3 in Fig.~\ref{delta2}. Again note that
there can be no element $a_0$ of $w$ between $L$ and $M$, or else
$[H,K,L,a_0,M,a_1,a_2,a_3]\sim P_3$. If $M$ participates in a pattern
$P_i$, then it forces $a_1,a_2,a_3$ also to participate. Since the
``1'' in $P_i$ cannot come from between $L$ and $M$, we conclude that
$L$ does not participate in $P_i$. But then we can replace $M$ by $L$
and argue that there is a pattern $P_i$ in $\widetilde{w}=d_M(w)$, a
contradiction.  Hence $M$ is pattern--free. \qed

\section{Relations among $\alpha,\,\,\beta,\,\,\gamma,\,\,\delta$
and $\epsilon$}
\label{recurrences}

Lemma~\ref{deletion1} shows that $h_n(x,k,l,m)$ does not depend on
$l$ or $m$, but rather on the differences $n-k$ and $x-k$. We shall
see further that there are only 5 ranges for $x-k$, that completely
determine the values of $h$: each of $x-k=0,1,2,3$ and $x-k\geq 4$
corresponds to exactly one of the 5 sequences listed earlier, and
defined as follows:
\[\left\{\begin{array}{l}
  h_n(x,x-0,l,m)=h_{n-x+0}(0,0,0,0)=:\alpha_{n-x-1}\\
  h_n(x,x-1,l,m)=h_{n-x+1}(1,0,0,0)=:\beta_{n-x+0}\\
  h_n(x,x-2,l,m)=h_{n-x+2}(2,0,0,0)=:\gamma_{n-x+1}\\
  h_n(x,x-3,l,m)=h_{n-x+3}(3,0,0,0)=:\delta_{n-x+2}\\
  h_n(x,x-\bar{x},l,m)=h_{n-x+\bar{x}}(\bar{x},0,0,0)=:
\epsilon_{n-x+3}
\end{array}\right\}
\Leftrightarrow \left\{\begin{array}{l} \alpha_n:=
h_{n+1\phantom{-\bar{x}}}(0,0,0,0)\\
\beta_{n}:=h_{n+1\phantom{-\bar{x}}}(1,0,0,0)\\ 
\gamma_n:=h_{n+1\phantom{-\bar{x}}}(2,0,0,0)\\
\delta_n:=h_{n+1\phantom{-\bar{x}}}(3,0,0,0)\\
\epsilon_n:=h_{n+\bar{x}-3}(\bar{x},0,0,0)
\end{array}\right\}\]
where $\bar{x}\geq 4$.  All of these definitions are justified by
Lemma~\ref{deletion1}, except for the definition of $\epsilon_n$,
which depends only on $n$ but not on $\bar{x}$ and whose explanation
will be given later. In our search for relations between these
sequences, it will be useful to define each sequence via an
alternative description in the terminology of the families
$\mathcal{F}$ discussed in Section~\ref{generating}.

\subsection{The sequence $\alpha$.} 
\begin{figure}[h]
$$\psdraw{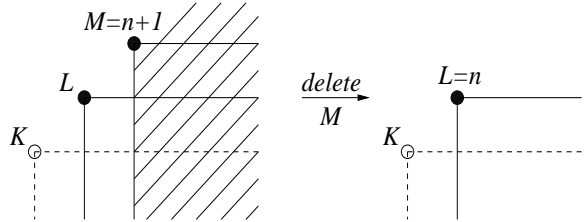}{3in}{1.15in}$$
\caption{Alternative Description of $\alpha_n$}
\label{alpha}
\end{figure}

By definition,
$\alpha_{n}=h_{n+1}(0,0,0,0)$. $M=n+1$ is the largest and the
rightmost element of any $w\in \mathcal{H}_{n+1}(0,0,0,0)$, and hence
it is pattern--free by (P2) (see Fig.~\ref{alpha}.)  Deleting $M$
results in a family of $321$--hex permutations of length $n$ without
any further restrictions, i.e. we have a bijection
\[d_M:\displaystyle{\mathcal{H}_{n+1}(0,0,0,0)\stackrel{\sim}{\rightarrow}
\mathcal{H}_n=\!\!\bigcup_{x,k,l,m}\!\!\mathcal{H}_n(x,k,l,m)}.\] This
justifies the alternative description of $\alpha_n$, stated in the
Introduction:
\begin{equation}
\alpha_n\,=\#\,\big\{\text{321--hex permutations
in\,\,}S_n\big\}= |\mathcal{H}_n|.
\end{equation}

\subsection{The sequence $\beta$.} 
\begin{figure}[h]
$$\psdraw{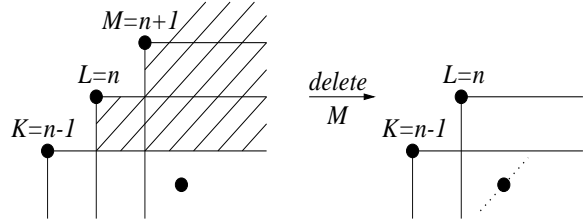}{3in}{1.15in}$$
\caption{Alternative Description of $\beta_{n}$}
\label{beta}
\end{figure}

By definition,
$\beta_{n}=h_{n+1}(1,0,0,0)$. Here $M=n+1$ is second from right to
left in $w\in \mathcal{H}_{n+1}(1,0,0,0)$, and by (P3) it is
pattern--free (see Fig.~\ref{beta}.) Thus, deleting $M$ imposes only
one extra condition: in the new 321--hex permutation
$\widetilde{w}\in\mathcal{H}_n$, the rightmost element in
$\mathcal{B}_1(\widetilde{w})$ is a ${K\backslash L}$--element or
lower, i.e. the largest two elements (the original $L=n$ and $K=n-1$)
belong to $\mathcal{B}_2(\widetilde{w})$.  This justifies the
alternative description:
\begin{equation}
\beta_{n}=\#\big\{w\in\mathcal{H}_n\,\big|\,
\{n,n-1\}\subset \mathcal{B}_2({w})\big\}\,\,\,
\text{for}\,\,n\geq 3.
\label{beta-def}
\end{equation}

In order to find a recursive description of $\beta_{n}$, note that
each $w\in \mathcal{H}_n$ either ends in $n$ or $n-1$, or both $n$ and
$n-1$ belong to $\mathcal{B}_2(w)$:
\begin{equation}
\alpha_n=\underbrace{\#\{w\in\mathcal{H}_n\,\big|\,
n\,\,\text{last}\}}_{\displaystyle{\alpha_{n-1}}}+
\underbrace{\#\{w\in\mathcal{H}_n\,\big|\,
n-1\,\,\text{last}\}}_{\displaystyle{\alpha_{n-1}}}+ \beta_{n}.
\label{beta-recursive1}
\end{equation}
By (P2), if $n$ or $n-1$ is the last element in $w\in\mathcal{H}_n$,
then it is pattern--free, so deleting it results in a bijection:
\begin{equation*}
d_n:\{w\in\mathcal{H}_n\,\big|\,
n\,\,\text{last}\}\stackrel{\sim}{\rightarrow}
\mathcal{H}_{n-1}\,\,\,\,\text{and}\,\,\,\,
d_{n-1}:\{w\in\mathcal{H}_n\,\big|\,
n-1\,\,\text{last}\}\stackrel{\sim}{\rightarrow} \mathcal{H}_{n-1}.
\label{beta-recursive2}
\end{equation*}
This justifies the use of $\alpha_{n-1}=|\mathcal{H}_{n-1}|$ in
(\ref{beta-recursive1}). Consequently, we have the recursive relation
$\beta_{n}=\alpha_n-2\alpha_{n-1}$ for $n\geq 3$.  \qed

\subsection{The sequence $\gamma$.}
\begin{figure}[h]
$$\psdraw{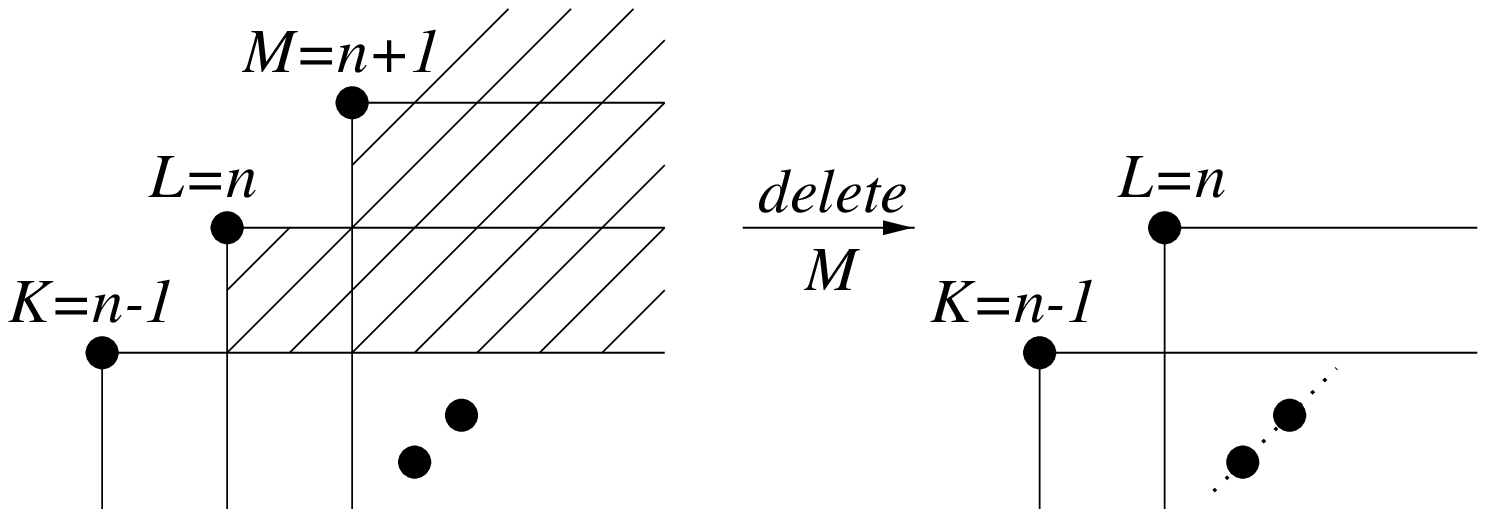}{3in}{1.15in}$$
\caption{$\gamma_{n}$}
\label{gamma}
\end{figure}

By definition, $\gamma_{n}=h_{n+1}(2,0,0,0)$. Here $M=n+1$ is third
from right to left in $w\in \mathcal{H}_{n+1}(2,0,0,0)$, and by (P3)
it is pattern--free (see Fig.~\ref{gamma}.) Deleting $M$ imposes the
following extra conditions: in the image
$\widetilde{w}\in\mathcal{H}_n$, the first two elements in
$\mathcal{B}_1(\widetilde{w})$ are ${K\backslash L}$--elements or
lower, i.e. the largest two elements (the original $L=n$ and $K=n-1$)
belong to $\mathcal{B}_2(w)$ and there are at least 2 numbers {\it
after} $L=n$.  This justifies the alternative description for $n\geq
4$:
\begin{equation}
\gamma_{n}=\#\big\{w\in\mathcal{H}_n\,\big|\,
\{w_s=n,\,\,w_t=n-1\}\subset \mathcal{B}_2(w),
\,\,t<s\leq n-2\big\}.
\label{gamma-def}
\end{equation}

\begin{figure}[h]
$$\psdraw{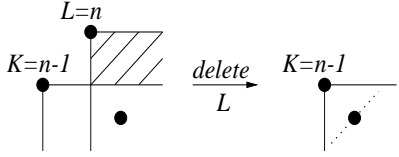}{2.05in}{0.8in}$$
\caption{Calculation of $\gamma_n$}
\label{gamma2}
\end{figure}

In order to find a recursive description of $\gamma_{n}$, note that
each $w$ counted in $\beta_{n}$ (i.e.
$\{n,\,\,n-1\}\subset\mathcal{B}_2(w)$), falls into one of the
following subcases: there is exactly one element {\it after} $n$
(hence in $\mathcal{B}_1(w)$), or there are at least two elements {\it
after} $n$ (hence in $\mathcal{B}_1(w)$):

\begin{equation}
\beta_{n}=\underbrace{\#\big\{w\in\mathcal{H}_n\,\big|\,
\{w_{n-1}=n,n-1\}\subset \mathcal{B}_2(w)\big\}}_
{\displaystyle{\alpha_{n-1}-\alpha_{n-2}}}+ \gamma_{n}.
\label{gamma-recursive1}
\end{equation}
The underbraced set in (\ref{gamma-recursive1}) is depicted in the LHS
of Fig.~\ref{gamma2}. By (P3), deletion of $L=n$ results in the
numerically equivalent set $S$ in the RHS of Fig.~\ref{gamma2}.  The
permutations in $S$ can be described as having their largest element
$K=n-1\in\mathcal{B}_{2}(\widetilde{w})$.  On the other hand,
$\mathcal{H}_{n-1}$ breaks into two disjoint groups: group $A$
consists of the permutations having the largest element $n-1$ in last
position, and group $B$ is the set $S$: $\mathcal{H}_{n-1}=A\sqcup S$.
Finally, note that $d_{n-1}:A\stackrel{\sim}{\rightarrow}
\mathcal{H}_{n-2}$, so that $|S|=\alpha_{n-1}-\alpha_{n-2}$. This
justifies the underbrace notation in (\ref{gamma-recursive1}), and
implies the formulas:
$\gamma_{n}=\beta_{n}-(\alpha_{n-1}-\alpha_{n-2})= \alpha_{n} -
3\alpha_{n-1} + \alpha_{n-2}$ for $n\geq 4$. \qed

\subsection{The sequence $\epsilon$.}
\begin{figure}[h]
$$\psdraw{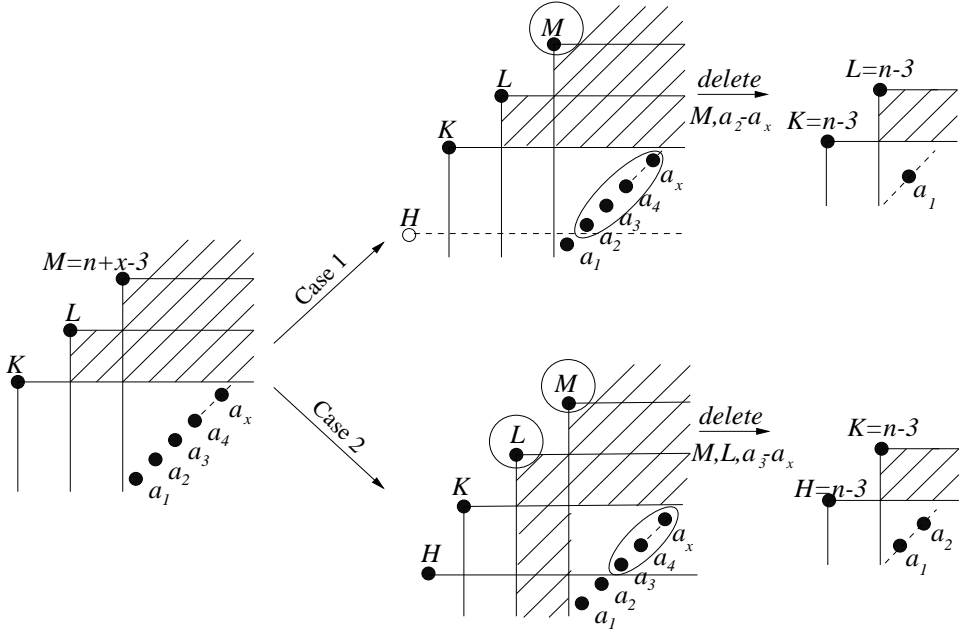}{5in}{3.3in}$$
\caption{Calculation of $\epsilon_{n}$}
\label{epsilon2}
\end{figure}

By definition, $\epsilon_{n}=h_{n+x-3}(x,0,0,0)$ for $x\geq 4$.  Let
$a_1,...,a_x$ be the ${X\backslash K}$--elements of
$w\in\mathcal{H}_{n+x-3}(x,0,0,0)$, where $x\geq 4$. There are two
cases to consider (see Fig.~\ref{epsilon2}.)

\smallskip
\noindent{\bf Case 1.} Except for $M,L,K$, all other elements of
$\mathcal{B}_2(w)$ are smaller than $a_2$. Our drawing shows the next
element $H\in\mathcal{B}_2(w)$ s.t. $H<a_2$ ($H$ may not exist in
$w$.)  By (P4), $a_2,a_3,a_4,...,a_x$ and $M$ are pattern--free.
After deletion, the remaining configuration in $\mathcal{H}_{n-3}$ is
identical to the alternative description of $\beta_{n-3}$.

\smallskip
\noindent{\bf Case 2.} After $M,L,K$, the fourth element $H$ of
$\mathcal{B}_2(w)$ is between $a_2$ and $a_3$: $a_2<H<a_3$.  Recall
from (P5) that there can be no element $a_0$ of $w$ between the
vertical lines of $L$ and $M$; or else, $a_0<a_1$, and hence
$[H,K,L,a_0,M,a_1,a_2,a_3]\sim P_1$.  Further, $a_3,a_4,...,a_x$, $M$
and $L$, are pattern--free. After deletion, the remaining
configuration in $\mathcal{H}_{n-3}$ is identical to the alternative
description of $\gamma_{n-3}$.

\smallskip
Note that the case $H>a_3$ is not allowable, or else
$[H,K,L,M,a_1,a_2,a_3, a_4]\sim P_2$ or $\sim P_4$ (when $a_3<H<a_4$
or $H>a_4$, respectively.)

\smallskip
Incidentally, the above discussion shows that $h_{n+x-3}(x,0,0,0)$
does not depend on $x$ as long as $x\geq 4$, and hence justifies the
definition of $\epsilon_n$.  We conclude that
$\epsilon_{n}=\beta_{n-3}+\gamma_{n-3}=
2\alpha_{n-3}-5\alpha_{n-4}+\alpha_{n-5}$ for $n\geq 7$. \qed

\subsection{The sequence $\delta$.}
\begin{figure}[h]
$$\psdraw{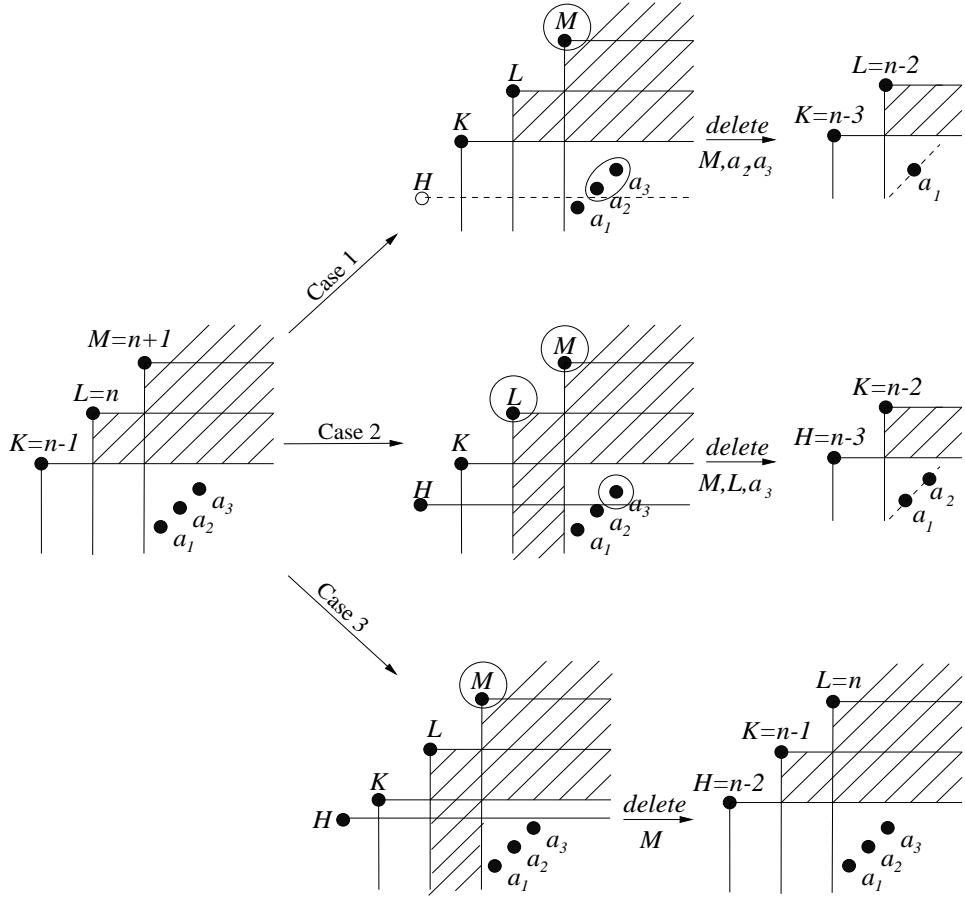}{5in}{4.7in}$$
\caption{Calculation of $\delta_{n}$}
\label{delta2}
\end{figure}

By definition, $\delta_{n}=h_{n+1}(3,0,0,0)$. As above, denote by
$a_1,a_2,a_3$ the ${X\backslash K}$--elements of
$w\in\mathcal{H}_{n+1}(3,0,0,0)$. There are three cases to consider
(see Fig.~\ref{delta2}.)

\noindent{\bf Case 1.} Except for $M,L,K$, all other elements of
$\mathcal{B}_2(w)$ are smaller than $a_2$. Our drawing shows the next
element $H\in\mathcal{B}_2(w)$ s.t. $H<a_2$ ($H$ may not exist in
$w$.) By (P4), $a_2,a_3$ and $M$ are pattern--free.  The remaining
configuration in $\mathcal{H}_{n-2}$ is identical to the alternative
description of $\beta_{n-2}$.

\smallskip
\noindent{\bf Case 2.} After $M,L,K$, the fourth element $H$ of
$\mathcal{B}_2(w)$ is between $a_2$ and $a_3$: $a_2<H<a_3$. By (P5),
there can be no element $a_0$ of $w$ between the vertical lines of $L$
and $M$; or else, $a_0<a_1$, $a_0\in \mathcal{B}_2(w)$, and hence
$[H,K,L,a_0,M,a_1,a_2,a_3]\sim P_1$. Further, $a_3$, $M$ and $L$ are
pattern--free. After deletion, the remaining configuration in
$\mathcal{H}_{n-2}$ is identical to the alternative description of
$\gamma_{n-2}$.

\smallskip
\noindent{\bf Case 3.} In contrast to the discussion of $\epsilon$, in
the case of $\delta$ it is possible to have $H>a_3$ as long as there
is no element $a_0$ of $w$ between the vertical lines of $L$ and $M$;
otherwise, $[H,K,L,a_0,M,a_1,a_2,a_3]~\sim~P_3$. By (P6), the largest
element $M=n+1$ is pattern--free.  After its deletion, the remaining
configuration in $\mathcal{H}_{n-1}$ is identical to the original
description of $\delta_{n-1}$.

\smallskip
We conclude that $\delta_{n}=\beta_{n-2}+\gamma_{n-2}+\delta_{n-1}=
\epsilon_{n+1}+\delta_{n-1}$ for $n\geq 6$. \qed

\medskip
We summarize the results in this Section in

\begin{lem} The sequences 
$\alpha,\,\,\beta,\,\,\gamma,\,\,\delta$ and $\epsilon$ satisfy the
following recursive relations:
\[\left|\begin{array}{ll}
\beta_n = \alpha_{n}-2\alpha_{n-1} & \text{for}\,\, n\geq 3; \\
\gamma_n = \alpha_{n} - 3\alpha_{n-1} + \alpha_{n-2} & 
\text{for}\,\, n\geq 4;\\
\delta_n = \epsilon_{n+1} + \delta_{n-1} & \text{for}\,\, n\geq 5;\\
\epsilon_n = 2\alpha_{n-3}-5\alpha_{n-4}+\alpha_{n-5} & 
\text{for}\,\, n\geq 6.
\end{array}\right.\]
\end{lem}

\section{Enumeration of 321--hex permutations}
\label{results}

We are now in a position to combine the recurrence formulas for
$\alpha,\beta,\gamma,\delta$ and $\epsilon$ into a single recurrence
for $\alpha$.  We first use the interpretation of $\alpha_n$ as
$|{\mathcal H}_n|$ and expand this in terms of the individual values
of $h_n(x,k,l,m)$. Thus, for a fixed $n$
\begin{eqnarray*}
\alpha_n&=& \sum_{x,k,l,m} h_n(x,k,l,m),
\end{eqnarray*}
where the sum is taken over $n\geq x\geq k\geq l\geq m\geq 0,\,\,x\not
= 0$.  Next, we break the sum into five separate sums depending on the
value of $x-k$; each such sum corresponds to the definitions of
$\alpha,\beta,\gamma,\delta$ and $\epsilon$, respectively. Note that
the sum for $\alpha$ (where $x=k$) requires two extra special cases
for $x=k=n-1$ and $x=k=n$; in both cases $h_n(x,k,l,m)=1$.
\begin{eqnarray*}
\Rightarrow\,\, \alpha_n&=&
1+(n-1)+\sum_{\begin{array}{c} 
      \scriptscriptstyle{n-2\geq x>0}\\\scriptscriptstyle{x\geq l\geq m\geq 0}
                   \end{array}}h_n(x,x,l,m)+
        \hspace*{-7mm}\sum_{\begin{array}{c} 
      \scriptscriptstyle{n\geq x}\\\scriptscriptstyle{x-1\geq l\geq m\geq 0}
             \end{array}}\hspace*{-7mm}h_n(x,x-1,l,m)\\
        &+&\hspace*{-7mm}\sum_{\begin{array}{c} 
      \scriptscriptstyle{n\geq x}\\\scriptscriptstyle{x-2\geq l\geq m\geq 0}
             \end{array}}\hspace*{-7mm}h_n(x,x-2,l,m)
        +\hspace*{-7mm}\sum_{\begin{array}{c} 
      \scriptscriptstyle{n\geq x}\\\scriptscriptstyle{x-3\geq l\geq m\geq 0}
             \end{array}}\hspace*{-7mm}h_n(x,x-3,l,m)
        +\hspace*{-7mm}\sum_{\begin{array}{c} 
      \scriptscriptstyle{n\geq x}\\
           \scriptscriptstyle{x-4\geq k\geq l\geq m\geq 0}
             \end{array}}\hspace*{-7mm}h_n(x,k,l,m).\\
\end{eqnarray*}
In the next step, we replace the $h_n$'s by the appropriate values
of $\alpha,\beta,\gamma,\delta$ and $\epsilon$.  The case of $\alpha$
requires some care as we must observe the condition that not all
the numerical parameters $x,k,l,m$ can be simultaneously equal; indeed,
this only happens in the special case $x=k=l=m=n$, which was broken
out from the main sum earlier.
\begin{equation*}
\Rightarrow\,\, \alpha_n = n+\hspace*{-5mm}\sum_{\begin{array}{c} 
                   \scriptscriptstyle{n-2\geq x,x>m}\\
                   \scriptscriptstyle{x\geq l\geq m\geq 0}
                   \end{array}}\hspace*{-6mm}
                   \alpha_{n-x-1}+\hspace*{-5mm}\sum_{\begin{array}{c} 
                   \scriptscriptstyle{n\geq x}\\\scriptscriptstyle{x-1\geq l\geq m\geq 0}
             \end{array}}\hspace*{-6mm}\beta_{n-x}
        +\hspace*{-5mm}\sum_{\begin{array}{c} 
     \scriptscriptstyle{n\geq x}\\\scriptscriptstyle{x-2\geq l\geq m\geq 0}
             \end{array}}\hspace*{-6mm}\gamma_{n-x+1}
        +\hspace*{-5mm}\sum_{\begin{array}{c} 
          \scriptscriptstyle{n\geq x}\\\scriptscriptstyle{x-3\geq l\geq m\geq 0}
             \end{array}}\hspace*{-6mm}\delta_{n-x+2}
        +\hspace*{-5mm}\sum_{\begin{array}{c} 
                   \scriptscriptstyle{n\geq x}\\\scriptscriptstyle
{x-4\geq k\geq l\geq m\geq 0}
             \end{array}}\hspace*{-6mm}\epsilon_{n-x+3}.
\end{equation*}
We replace the indices $k$, $l$ and $m$ by binomial coefficients:
\begin{equation*}
\alpha_n=n+\sum_{i=1}^{n-1}\left({\scriptstyle\binom{n+1-i}{2}-1}\right) \alpha_i+
          \sum_{i=3}^{n-1}{\scriptstyle\binom{n+1-i}{2}}\beta_i
        +\sum_{i=4}^{n-1}{\scriptstyle\binom{n+1-i}{2}}\gamma_i
        +\sum_{i=5}^{n-1}{\scriptstyle\binom{n+1-i}{2}}\delta_i
        +\sum_{i=6}^{n-1}{\scriptstyle\binom{n+2-i}{3}}\epsilon_i.
\end{equation*}
\noindent Finally, we use the recurrences from Sect.~\ref{recurrences} to
obtain a summation in terms of $\alpha$ alone:
\begin{eqnarray*}
\alpha_n\!&=& n+
\sum_{i=1}^{n-1}\left({\scriptstyle\binom{n+1-i}{2}-1}\right)
\alpha_i+
\!\!\sum_{i=3}^{n-1}{\scriptstyle\binom{n+1-i}{2}}(\alpha_i-2\alpha_{i-1})
+\!\!\sum_{i=4}^{n-1}{\scriptstyle\binom{n+1-i}{2}}(\alpha_i-3\alpha_{i-1}+\alpha_{i-2})\\&&\phantom{n}+
\sum_{i=5}^{n-1}{\scriptstyle\binom{n+1-i}{2}}
\big(2\alpha_{i-2}-3\alpha_{i-3}
-2\sum_{j=3}^{i-4}\alpha_j-4\alpha_2+\alpha_1\big)+{\scriptstyle\binom{n-4}{2}}\delta_5\\&&\phantom{n}+
\sum_{i=6}^{n-1}
{\scriptstyle\binom{n+2-i}{3}}(2\alpha_{i-3}-5\alpha_{i-4}+\alpha_{i-5}).\\
\end{eqnarray*} 
We simplify this large expression into the following full--history,
linear recurrence relation with cubic polynomial coefficients, valid
for $n\geq 6$:
\begin{eqnarray*}
\alpha_n&=& 2\alpha_{n-1}+3\alpha_{n-2}+\sum_{k=3}^{n-4}Q(k)\alpha_{n-k}+R(n),
\end{eqnarray*}
where $Q(k)=-\frac{2}{3}k^3 + \frac{11}{2}k^2 -\frac{71}{6} + 9$ and
$R(n)=-\frac{17}{3}n^3+85n^2-\frac{1267}{3}n+704.$ The polynomial
$R(n)$ arises from three special cases for the polynomial coefficients
of $\alpha_3,\alpha_2$ and $\alpha_1$.  Since $\deg Q(n)=\deg R(n)=3$,
we need 4 successive history eliminations of the form
$\alpha_n-\alpha_{n-1}$; combined with the 3 initial terms
$\alpha_n,\alpha_{n-1}$ and $\alpha_{n-2}$, this produces the desired
{order--six constant--coefficient linear recurrence} for all $n\geq
6$:
\begin{equation}
\alpha_{n}=6\alpha_{n-1}-11\alpha_{n-2}+
9\alpha_{n-3}-4\alpha_{n-4}-4\alpha_{n-5}+\alpha_{n-6}.
\label{main-recurrence}
\end{equation}

\noindent Consequently, the number of the 321--hexagon avoiding
permutations of length $n$ is given by the formula:
\begin{equation}
c_1R_1^n+c_2R_2^n+c_3R_3^n+c_4R_4^n+c_5R_5^n+
\overline{c_5{R_5}^n},
\label{formula}
\end{equation}
where the roots and coefficients are listed in the Introduction. This
completes the proof of Theorem~\ref{main-result}.  \qed

\smallskip
The sixth--degree characteristic polynomial of our recurrence relation
(\ref{main-recurrence}) is irreducible over $\mathbb{Q}$ and has
Galois group $S_6$, as calculated by Maple.  This means that there are
no further algebraic relations among the roots $R_i$ in
(\ref{formula}), and thus we cannot hope for any better closed--form
results.

\smallskip
However, our numerical approximations of the roots and coefficients
can yield exact values for the number of 321--hex permutations for any
fixed length $n$, exploiting the fact that the value being
approximated is known to be an integer.  Furthermore, since the two
roots of modulus less than $1$ make such small contributions, they can
be dropped and the following formula is exact for all $n$:
\begin{equation}
\alpha_n = [c_3R_3^n+c_4R_4^n+c_5R_5^n+\overline{c_5R_5^n}],
\end{equation}
where the braces denote rounding to the nearest integer.

\section{Extensions and Further Discussion}
\label{extensions}

\begin{figure}[h]
$$\psdraw{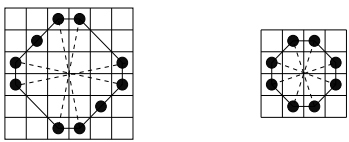}{1.8in}{0.7in}$$
\caption{$6\times 6$ and $4\times 4$ Patterns}
\label{extension}
\end{figure}

The study of the octagonal patterns $P_i$ in the present paper was
motivated by their apperance in the representation theory of $S_n$ via
heap--avoidance and thus Kazhdan--Lusztig polynomials and Schubert
varieties. We refer to the enumeration of $S_n(321,\mathcal{P})$ as
the ``$8\times 8$ case''. From a purely combinatorial viewpoint, it is
natural to ask what happens in the analogous smaller $6\times 6$
and $4\times 4$ cases whose generalized patterns are depicted in
Fig.~\ref{extension}. 

\noindent 
To obtain these cases, for each $i=1,2,3,4$ shorten by one the lengths
of both $\mathcal{B}_1(P_i)$ and $\mathcal{B}_2(P_i)$ by removing a
chosen fixed point appearing in all $P_i$'s. Thus, define
\begin{eqnarray*}
\mathcal{P}^{6}&=&\{[351624],[356124],[451623],[456123]\}\,\,\text{and}\\
\mathcal{P}^{4}&=&\{[2143],[3142],[2413],[3412]\}
\end{eqnarray*}
to be the families of avoided octagonal patterns in the $6\times 6$
and the $4\times 4$ cases, respectively. Both of these cases lead
again to linear recursive relations with constant coefficients.  The
proofs below follow closely the method described in the $8\times 8$
case, so we leave the details for verification to the reader.

\subsection{The $6\times 6$--case}

\begin{thm}
Let $\mathcal{H}_n^{6}=S_n(321,\mathcal{P}^{6})$, and
$\alpha_n=|\mathcal{H}_n^{6}|$. Then $\{\alpha_n\}$ satisfies a
6--term linear recursive relation with constant coefficients:
\begin{equation}
\alpha_{n+1}=4\alpha_n-4\alpha_{n-1}+3\alpha_{n-2}+\alpha_{n-3}-\alpha_{n-4}\,\,\,\,\text{for all}\,\,
n\geq 1.
\end{equation} Consequently, for all $n\geq 1$,
\begin{equation}
\alpha_n=c_1R_1^n+c_2R_2^n+c_3R_3^n+c_4R_4^n+\overline{c_4R_4^n}
\end{equation} where the roots and coefficients are rounded off below to 5
digits after the decimal point:
\[\begin{array}{lll}
R_1\approx -0.49569        & &c_1\approx  \phantom{-}0.63205\\
R_2\approx\phantom{-}0.51154& &c_2\approx \phantom{-}0.53110 \\
R_3\approx\phantom{-}3.03090& &c_3\approx \phantom{-}0.50154\\
R_4\approx\phantom{-}0.47662-1.03635i& &c_4\approx  -0.19482+0.11092i \\
\end{array}\]
\label{6-by-6}
\end{thm}
\vspace*{-4mm}
\noindent The degree--5 characteristic polynomial of the recurrence is
irreducible over $\mathbb{Q}$ and has the largest possible Galois
group, $S_5$, as calculated by Maple.  We can again drop the small
roots $R_1$ and $R_2$, rounding off the remainder to the nearest
integer:
\begin{equation}
\alpha_n=[c_3R_3^n+c_4R_4^n+\overline{c_4R_4^n}\big]\,\,\text{for all}\,\,
n\geq 1.
\end{equation}
The first values of $\alpha_n$ are:
1, 2, 5, 14, 42, 128, 389, 1179, 3572, 10825, 32810, 99446.

\medskip
\noindent{\sc Proof}(Theorem~\ref{6-by-6}): We modify the discussion
in the proof for the $8\times 8$ case. Define the generating function
$h_n(x,l,m)$ with one fewer parameter, thus taking into account only
the elements $M$ and $L$ of $w\in \mathcal{H}_n^{6}$. As in
Lemma~\ref{deletion1}, the $L$--elements are pattern--free, and
therefore $h_n(x,l,m)=h_{n-l}(x-l,0,0)$.

\smallskip
Next, define the sequences $\alpha_n=h_{n+1}(0,0,0)$,
$\beta_{n+1}=h_n(1,0,0)$, $\gamma_n=h_{n+1}(2,0,0)$ and
$\delta_n=h_{n+x-2}(x,0,0)$ for $x\geq 3$. Note that
$\alpha_n=h_{n+1}(0,0,0)=|\mathcal{H}_n^{6}|$ since $\alpha_n$ is the
number of $w\in \mathcal{H}_{n+1}^{6}$ in which the largest element
$n+1$ is at the end and hence it is pattern--free. Further, the relations
among the sequences are:
\[\left|\begin{array}{ll}
\beta_n=\alpha_n\phantom{\scriptscriptstyle{-1}}-
\alpha_{n-1}&\text{for}\,\,n\geq 2;\\
\gamma_n=\alpha_{n-2}+\beta_{n-2}+\gamma_{n-1}&\text{for}\,\,n\geq
3;\\ \delta_n=\alpha_{n-3}+\beta_{n-3}&\text{for}\,\,n\geq 4.
\end{array}\right.\]
Now we are ready to express everything in terms of $\alpha_n=
\sum_{x,l,m} h_n(x,l,m)$:
\[\alpha_n=1+\sum h_n(x,x,m)+\sum h_n(x,x-1,m)+
\sum h_n(x,x-2,m)+\sum h_n(x,l,m).\]
Here the ``1'' counts the identity permutation: no $M$--element in $w$;
in the first sum $L$ is smaller than all $X$--elements (or $L$ does
not exist); in the second sum $L$ is larger than exactly one
$X$--element; in the third sum $L$ is larger than exactly two
$X$--elements; and in the forth sum $L$ is larger than three or more
$X$--elements and hence $x-3\geq l$.
\begin{eqnarray*}
&\Rightarrow&
\alpha_n=1+\sum_{l=1}^{n-1}l\alpha_{n-l-1}+\sum_{l=0}^{n-3}
(l+1)\beta_{n-l-1}+\sum_{l=0}^{n-4}(l+1)\gamma_{n-l-1}+
\sum_{l=0}^{n-5}{\scriptscriptstyle\binom{l+2}{m}}\delta_{n-l-1}\\
&\Rightarrow&
\alpha_{n+1}-\alpha_{n}=\sum_{l=1}^{n}\alpha_{n-l}+
\sum_{l=0}^{n-2}\beta_{n-l}+\sum_{l=0}^{n-3}\gamma_{n-l}+
\sum_{l=0}^{n-4}(l+1)\delta_{n-l}\\
&\Rightarrow&
\alpha_{n+2}-2\alpha_{n+1}+\alpha_{n}=\alpha_n+\beta_{n+1}+\gamma_{n+1}+
\sum_{l=0}^{n-3}\delta_{n-l+1}\\
&\Rightarrow&
\alpha_{n+3}-3\alpha_{n+2}+2\alpha_{n+1}=(\beta_{n+2}-\beta_{n+1})
+(\gamma_{n+2}-\gamma_{n+1})+\delta_{n+2}.
\end{eqnarray*}

\medskip
\noindent Conveniently, each summand on the RHS can be expressed in
terms of $\alpha$, including
$\gamma_{n+2}-\gamma_{n+1}=\alpha_n+\beta_n=2\alpha_n-\alpha_{n-1}$:
\[\Rightarrow\,\,\alpha_{n+3}=4\alpha_{n+2}-4\alpha_{n+1}+3\alpha_n+\alpha_{n-1}-\alpha_{n-2}.\qed\]

\subsection{The $4\times 4$ case.}

\begin{thm}
Let $\mathcal{H}_n^{4}=S_n(321,\mathcal{P}^{4})$, and
$\alpha_n=|\mathcal{H}_n^{4}|$. Then
\begin{equation}
\alpha_{n}=(n-1)^2+1\,\,\,\,\text{for all}\,\,n\geq 1.
\end{equation}
In particular, $\{\alpha_n\}$ satisfies the 4--term linear recursive relation
\begin{equation}
\alpha_{n+1}=3\alpha_n-3\alpha_{n-1}+\alpha_{n-2}\,\,\,\,
\text{for all}\,\,n\geq 1.
\end{equation}
\end{thm}
\noindent Note that in contrast to the $6\times 6$ and $8\times 8$
cases, the characteristic polynomial here factors (completely) over
$\mathbb{Q}$: $(x-1)^3$.

\smallskip
\noindent
{\sc Proof:} As expected, we define the generating function $h_n(x,m)$
with only two parameters, keeping track only of $M$ in $w\in
\mathcal{H}_n^{4}$. Obviously, the $M$--elements are pattern--free, so
$h_n(x,m)=h_{n-m}(x-m,0)$. Further, letting $\beta_n=h_{n+1}(1,0)$
and $\gamma_n=h_{n+x-1}(x,0)$ for $x\geq 2$, one quickly
discovers that $\beta_n=\beta_{n-1}+1$, so $\beta_n=n+1$ for $n\geq 3$,
and $\gamma_n=n-1$ for $n\geq 1$. Express $\alpha_n=
\sum_{x,m} h_n(x,m)$ as the sum
\[\alpha_n=1+\sum h_n(x,x-1)+\sum h_n(x,m).\]
As before, here the ``1'' counts the identity permutation: no
$M$--element in $w$; in the first sum $M$ is larger than exactly one
$X$--element; and in the third sum $M$ is larger than two or more
$X$--elements and hence $x-2\geq m$.
\begin{eqnarray*}
\Rightarrow&&
\alpha_n=1+\sum_{m=0}^{n-2}\beta_{n-m-1}+\sum_{m=0}^{n-3}\gamma_{n-m-1}\\
\Rightarrow&&\alpha_{n+1}-\alpha_n=\beta_n+\gamma_n\,\,
\Rightarrow\,\,
\alpha_{n+2}-2\alpha_{n+1}+\alpha_n=2\\
\Rightarrow&&
\alpha_{n+3}-3\alpha_{n+2}+3\alpha_{n+1}+\alpha_{n}=0\,\,\text{for all}\,\,n\geq 1. \qed
\end{eqnarray*}

\subsection{Further Discussion} 
Until now, there were relatively few known examples of sets of
permutations whose avoidance led to linear, polynomial, or
exponential formulas (see \cite{St1,We2}.) After the
successful enumeration of the $4\times 4$, $6\times 6$ and the
321--hex $8\times 8$ cases, it is tempting to generalize the
recursive sequence method of this paper to the corresponding larger
sets of $2k\times 2k$ patterns. At this point, it is not surprising to
conjecture that all these families yield linear recursive relations
with constant coefficients.  In fact, when the result of the present
paper was publicized, Herbert Wilf requested that many more examples
of such ``linear'' families be found. These and other related
questions will be answered positively in a forthcoming paper by
Stankova--Frenkel.

\section* {Acknowledgments}

The authors would like to thank Sara Billey and Gregory Warrington for
sharing their results on Kazhdan--Lusztig polynomials for
321--hexagon--avoiding permutations and suggesting the problem on
enumeration of forbidden subsequences in this paper.


\begin{thebibliography}{[14]}

\bibitem{BW} E. Babson, J. West, {The permutations $123p_4...p_t$
and $321p_4...p_t$ are Wilf--equivalent}, Graphs Comb 16 (2000) 4,
373--380.

\bibitem{BWX} J. Backelin, J. West and G. Xin, {Wilf--equivalence for
singleton classes}, in {Proc. 13th Conf. in Formal Power Series
and Algebraic Combinatorics}, Tempe 2001.

\bibitem{BJS} S. Billey, W. Jockusch, and R. Stanley,
{Some combinatorial properties of Schubert polynomials,}
J. Alg. Comb., 2 (1993), 345--374.

\bibitem{Billey} S. Billey and G. Warrington.  {Kazhdan--Lusztig
polynomials for 321--hexagon--avoiding permutations,}
arXiv:math.CO/0005052, 5 May 2000.

\bibitem{BS} R. Bott and H. Samelson, {Applications of the theory
of Morse to symmetric spaces}, Amer. J. Math., 80 (1958), 964--1029.

\bibitem{BK} J.-L. Brylinski and M. Kashiwara, {Kazhdan--Lusztig
Conjectures and Holonomic Systems}, Invent. Math., 64 (1981),
387--410.

\bibitem{CGHK} F.R.K. Chung, R.L. Graham, V.E. Hoggatt Jr. and
M. Kleiman, {The Number of Baxter Permutations}, J. Combin.
Theory Ser. A 24 (1978), 382-394.

\bibitem{Hum} J. Humphreys, {Reflection Groups and Coxeter
Groups}, Cambridge University Press, 1990.

\bibitem{KL} D. Kazhdan and G. Lusztig, {Representations and
Coxeter Groups and Hecke Algebras}, Inv. Math., 53 (1979), 165--184.

\bibitem{KL2} D. Kazhdan and G. Lusztig, {Schubert Varieties and
Poincar\'{e} Duality}, Proc. Symp. Pure. Math., AMS, 36 (1980),
185--203.

\bibitem{Schensted} C. Schensted, {Longest Increasing and
Decreasing Subsequences}, Canad. J. Math. 13 (1961), pp. 179--191.

\bibitem{St1} Z. Stankova,  {Forbidden subsequences},
Disc. Math. 132 (1994) 291--316.

\bibitem{St2} Z. Stankova,  {Classification of forbidden
subsequences of length 4}, Europ. J. Combinatorics (1996) 17, 501--517.

\bibitem{St3} Z. Stankova--Frenkel and J. West, {A New Class of
Wilf--Equivalent Permutations}, submitted to J. Alg. Comb.;
arXiv:math.CO/0103152 (2001).

\bibitem{Ti} J. Tits, {Le probl\`{e}me des mots dans les groupes
de Coxeter,} Symposia Math, 1 (1968), 175--185. Ist. Naz. Alta
Mat. (1968), {Symposia Math.,} Vol. 1, Academic Press, London.

\bibitem{Vi} G. Viennot, {Heaps of pieces. I. Basic definitions
and combinatorial lemmas}, Graph theory and its applications: East and
West, Jinan, (1986), 542--570.

\bibitem{We1} J. West,  {Generating trees and the Catalan and
Schr\"{o}der numbers}, Disc. Math. 146 (1995) 247--262.

\bibitem{We2} J. West,  {Generating trees and forbidden
subsequences}, Disc. Math. 157 (1996) 363--374.


\end{thebibliography}
\end{document}